\documentclass[12pt]{article}
\usepackage{amsmath,amssymb,hhline}
\usepackage[russian]{babel}

\usepackage[T2A]{fontenc}    
\usepackage[cp1251]{inputenc}
\usepackage[russian]{babel} 

\begin{document}

\begin{center}
{\large\bf R.~M.~Trigub}
\end{center}

\begin{center}
{\large \bf On the best approximation of constants\\ by polynomials with integer coefficients}
\end{center}

\emph{\textbf{Abstract.}} In this paper, exact rate of decrease of best approximations of non-integer numbers by polynomials with integer coefficients of the growing exponentials is found on a disk in complex plane, on a cube in $\mathbb{R}^d$, and on a ball in $\mathbb{R}^d$. While in the first two cases the $\sup$-norm is used, the third one is fulfilled in $L_p$, $1\leq p<\infty$.

Comments are also given (two remarks in the end of the paper).



\begin{center}
{\large\bf Р.~М.~Тригуб}
\end{center}

\begin{center}
{\large \bf О наилучшем приближении констант\\ многочленами с целыми коэффициентами}
\end{center}

\emph{\textbf{Аннотация.}} В статье найден точный порядок убываний наилучших приближений нецелых чисел многочленами с целыми коэффициентами растущих степеней на круге в комплексной плоскости, на кубе в $\mathbb{R}^d$ и на шаре в $\mathbb{R}^d$. В первых двух случаях используется $\sup$-норма, а в третьем --- $L_p$, $1\leq p<\infty$. Приведены комментарии (два замечания в конце статьи).


\emph{\textbf{Ключевые слова и фразы:}} трансфинитный диаметр, многочлен Чебышева, экстремальные свойства многочленов, неравенство Маркова, теоремы Бернштейна, наилучшее приближение, целые алгебраические числа, $q$-ичные дроби.

\begin{center}
    \textbf{\S1 Введение}
\end{center}

С.~Н.~Бернштейн на первом математическом съезде СССР (1930) поставил вопрос о наилучшем приближении нецелого числа многочленами с целыми коэффициентами растущей степени (см. [\ref{Bernsht}], I, 43 и 47). Сразу же Р.~О.~Кузьмин и Л.~В.~Канторович указали некоторую оценку приближения сверху, не зависящую от природы числа [\ref{Kantor}]. Полностью вопрос решен для отрезка $[\delta,1-\delta]$, $\displaystyle \delta\in\big(0,1/2\big)$, в статье [\ref{Trigub62}] (см. также [\ref{Trigub_Belinsky}], \textbf{5.4.16}). А в работе [\ref{Trigub2001}] приведены несколько случаев точного порядка убывания наилучших приближений констант на отрезке $[\alpha,\beta]$, $0<\alpha<\beta<1$, в зависимости от арифметической природы числа и отрезка. В этой же статье впервые рассмотрены вопросы о приближении гладких функций и констант многочленами с натуральными коэффициентами (на отрезке, лежащем на отрицательной полуоси $\mathbb{R}$).

Очевидно, что если какая-либо функция, отличная от многочлена, допускает равномерное приближение многочленами $q_n$ с целыми коэффициентами на компакте $K\subset \mathbb{C}$, напр., то существует многочлен $X$ с целыми коэффициентами такой, что
\begin{equation*}
  0<\max\limits_{z\in K}|X(z)|<1.
\end{equation*}

Следовательно, трансфинитный диаметр $K$ меньше единицы (см., напр., [\ref{Timan}], \textbf{2.13}, \textbf{13}). Если же существует такой многочлен с условием $0<|X(z)|<1$, то любая функция, которую можно приблизить многочленами $p_n$  с любыми коэффициентами допускает приближение многочленами $q_n$ (достаточно приблизить константу $\lambda=\frac{1}{2})$. Если же на компакте $K$ есть целочисленная точка (число из $\mathbb{Z}+i\mathbb{Z}$) или, более общо, целые алгебраические числа вместе со всеми сопряженными, то функция должна удовлетворять некоторым арифметическим условиям. Напр., вещественная непрерывная функция на $[-1,1]$ является пределом многочленов $q_n$ при $n\rightarrow\infty$ в том и только в том случае, если $f(0)$ и $ \frac{1}{2}\big(f(-1)\pm f(1)\big)\in\mathbb{Z}$ (см. [\ref{Lorentz}], гл.~2, \S4). Отметим, что трансфинитный диаметр отрезка равен четверти его длины, круга --- радиусу, эллипса --- полусумме полуосей.

В настоящей статье изучаются наилучшие приближения разных констант в случаях, когда компакт $K$ --- это круг в $\mathbb{C}$ (\S2, теорема~1 и следствие), куб в $\mathbb{R}^d$ (\S3, теорема~2) и шар в $\mathbb{R}^d$ (\S4, теорема~3). В случае шара --- интегральная метрика и арифметические условия отпадают.

Оценка приближения в пространстве $C$ снизу существенно зависит от ответа на вопрос о точном порядке роста многочленов $p_n$ при $n\rightarrow\infty$ вне компакта $K$ при данной $\sup$-норме на $K$. В случае $K=[a,b]\subset\mathbb{R}$, напр., как давно известно (см. [\ref{Timan}], \textbf{2.13}, \textbf{14}), для любого многочлена $p_n$ степени $n$ при $x\in\mathbb{R}\setminus[a,b]$
\begin{equation*}
  |p_n(x)|\leq|C_n(x)|\max\limits_{[a,b]}|p_n(x)|,
\end{equation*}
где $C_n(x)=C_n(x;a,b)$ --- многочлен Чебышева для $[a,b]\subset\mathbb{R}$:
\begin{equation*}
  C_n(x)=\frac{1}{2}\Bigg[\Bigg(\frac{2x-a-b}{b-a}+\sqrt{\Big(\frac{2x-a-b}{b-a}\Big)^2-1}\Bigg)^n+\Bigg(\frac{2x-a-b}{b-a}-\sqrt{\Big(\frac{2x-a-b}{b-a}\Big)^2-1}\Bigg)^n\Bigg].
\end{equation*}

А если $a>0$, напр., то при некотором $\displaystyle\theta\in\Big[\frac{1}{2},1\Big]$, как следует из предыдущего,
\begin{equation}\label{eq_1}
  |C_n(0)|=\theta\Big(\frac{\sqrt{b}+\sqrt{a}}{\sqrt{b}-\sqrt{a}}\Big)^n.
\end{equation}

Здесь и далее $p_n$ --- многочлен степени не выше $n$ с любыми коэффициентами, а $q_n$ --- только с целыми. Через $c(\alpha,\beta)$ с разными индексами будем обозначить некоторые положительные величины, зависящие только от $\alpha$ и $\beta$. См. также комментарии 1 и 2 в конце статьи.

\begin{center}
    \textbf{\S2 Приближение констант на круге в $\mathbb{C}$}
\end{center}

Пусть $\lambda\in\mathbb{C}$ и не является целым числом, т.е. $\lambda\notin\mathbb{Z}+i\mathbb{Z}$. В круге
\begin{equation*}
  K_r=K_r(z_0)=\big\{z\in\mathbb{C}: \ |z-z_0|\leq r\big\}
\end{equation*}
нет целочисленных точек и, значит, $\displaystyle 0<r<\frac{1}{\sqrt{2}}$. После сдвига на целое число можно считать, что $z_0\in\Pi_1$, где $\Pi_1$ --- замкнутый квадрат с вершинами 0, 1, $1+i$, $i$. Если $\displaystyle \textrm{Re}~z_0>\frac{1}{2}$, то после замены $z$ на $1-z$ получим $\displaystyle \textrm{Re}~z_0<\frac{1}{2}$. По той же причине можно считать далее, что и $\displaystyle \textrm{Im}~z\leq\frac{1}{2}$, т.е. центр круга $\displaystyle z_0\in\frac{1}{2}\Pi_1$ и $0$ является ближайшей к $K_r$ вершиной $\Pi_1$.

\medskip
\textbf{Теорема 1} \emph{При сделанных предположениях и любом $n$}
\begin{equation*}
  E_n^e(\lambda;K_r)=\min\limits_{q_n}\max\limits_{K_r}|\lambda-q_n(z)|\leq(n+1)\rho^n,
\end{equation*}
\emph{где $\rho=\max\{\rho_1,\rho_2\}$, а $\displaystyle\rho_1=\frac{r}{|z_0|}$ и $\rho_2=|z_0|+r$.}

\emph{ При этом $\rho$ уменьшить нельзя, в общем случае.}

$\triangleright$ Так как
\begin{equation*}
  \Big(\frac{z-z_0}{-z_0}\Big)^m=1+\sum\limits_{k=1}^m a_kz^k,
\end{equation*}
то
\begin{equation*}
  \max\limits_{K_r}\Big|\lambda-\sum\limits_{k=1}^n\lambda a_kz^k\Big|\leq|\lambda|\Big(\frac{r}{|z_0|}\Big)^n=|\lambda|\rho_1^n.
\end{equation*}

Заменяя $\lambda a_1$ ближайшим целым числом $c_1$ и применяя к их разности ("дробная часть") это же неравенство с заменой $n$ на $n-1$, получаем
\begin{equation*}
    \max\limits_{K_r}\Big|\lambda-c_1z-\sum\limits_{k=2}^n b_kz^k\Big|\leq |\lambda|\rho_1^n+\Big(\frac{r}{|z_0|}\Big)^{n-1}\max\limits_{K_r}|z|=|\lambda|\rho_1^n+\rho_1^{n-1}\rho_2.
  \end{equation*}

Продолжая таким же образом, получаем, считая $|\lambda|\leq1$,
\begin{equation*}
  E_n^e(\lambda;K_r)\leq\max\limits_{K_r}|\lambda-q_n(z)|\leq\sum\limits_{k=0}^n\rho_1^{n-k}\rho_2^k\leq(n+1)\max\big\{\rho_1^n,\rho_2^n\big\}=(n+1)\rho^n.
\end{equation*}

Очевидно, что при $\rho_1\neq\rho_2$ множитель $(n+1)$ можно заменить ограниченной по $n$ величиной.

А теперь для оценки приближения снизу воспользуемся неравенством: для любого многочлена $p_n$ при $|z-z_0|>r$
\begin{equation*}
  |p_n(z)|\leq\Big(\frac{|z-z_0|}{r}\Big)^n\max\limits_{K_r}|p_n(z)|
\end{equation*}
(см. [\ref{Timan}], [\ref{Dzjad}], [\ref{Trigub_Belinsky}]).

При некотором многочлене $q_n$ имеем
\begin{equation*}
  E_n^e(\lambda;K_r)=\max\limits_{K_r}|\lambda-q_n(z)|\geq|\lambda-q_n(0)|\Big(\frac{r}{|z_0|}\Big)^n\geq c(\lambda)\rho_1^n.
\end{equation*}

Это оценка приближения снизу при $\rho_1\geq\rho_2$.

Пусть теперь $\rho_2>\rho_1$, т.е. $\displaystyle r<\frac{|z_0|^2}{1-|z_0|}$.

Предположим, что $\displaystyle|z_0|+r=\frac{1}{q}$, где $q\in \mathbb{N}$ и $q\geq2$, а $\displaystyle \frac{1}{q}e^{i\varphi}\in K_r$, где угол $\varphi$ такой, что круг $K_r$ пересекается с интервалом $(0,1)$ по некоторому отрезку, т.е. $\displaystyle 0\leq\sin\varphi<\frac{r}{|z_0|}$.

Возьмём на этом отрезке число $\displaystyle \frac{p}{q}$ ($p\in \mathbb{N}$, а далее $p_1$ --- целое). Тогда
\begin{equation*}
  E_n^e(\lambda;K_r)\geq\min\limits_{q_n}\Big|\lambda-q_n\Big(\frac{p}{q}\Big)\Big|\geq\min\limits_{p_1}\Big|\lambda-\frac{p_1}{q^n}\Big|\geq\frac{1}{q^{n+2}},
\end{equation*}
если в $q$-ичном разложении числа $\textrm{Re}~\lambda$ не следуют подряд два раза 0 или $q-1$.

\vskip 1mm
Но $\displaystyle \frac{1}{q}=|z_0|+r=\rho_2$. Так что есть числа $\lambda$ и круги $K_r$, при которых и $\rho_2$ уменьшить нельзя. \qquad\qquad\qquad\qquad\qquad\qquad\qquad\qquad\qquad\qquad\qquad\qquad\qquad\qquad\qquad\qquad\qquad$\blacktriangleleft$

\medskip
Рассмотрим ещё только особый случай $\displaystyle \lambda=\frac{p}{q^s}$ ($p$ --- целое, $s$ и $q\in\mathbb{N}$, $q\geq2$).

Положим
\begin{equation*}
  q_n(z)=\frac{1}{q}-\frac{1}{q}(1-qz)^n.
\end{equation*}

Тогда при $s=1$
\begin{equation*}
  \lambda-pq_n(z)=\frac{p}{q}-pq_n(z)=\lambda(1-qz)^n
\end{equation*}
и имеем оценку приближения $\lambda$.

А далее применяем индукцию по $s$, учитывая, что
\begin{equation*}
  \Big(\frac{1}{q}-q_n(z)\Big)^s=\frac{1}{q^s}(1-qz)^{sn}
\end{equation*}
и после умножения на $p$ и применения формулы бинома Ньютона
\begin{equation*}
  \frac{p}{q^s}-s\frac{p}{q^{s-1}}q_n(z)+...+(-1)^spq_n^s(z)=\frac{p}{q^s}(1-qz)^{sn}.
\end{equation*}

Заменяя $n$ на $\displaystyle \Big[\frac{n}{s}\Big]$, получаем
\begin{equation*}
  E_n^e\Big(\frac{p}{q^s};K_r\Big)\leq c(\lambda,s)\max\limits_{K_r}|1-qz|^n=c(\lambda,s)q^n\max\limits_{K_r}\Big|z-\frac{1}{q}\Big|^n.
\end{equation*}

А это может быть и меньше $\rho^n$.

Отметим теперь, что от приближения констант можно перейти к приближению любых функций, допускающих приближение многочленами $p_n$.

В примере с кругом $K_r(z_0)$ установим связь между наилучшими приближениями многочленами $p_n$ с любыми коэффициентами с наилучшими приближениями многочленами $q_n$.

Пусть функция $f$ аналитична в замкнутом круге $K_R(z_0)$ с тем же центром $z_0$ при некотором радиусе $R>r$, а $s$ --- наименьшее число из $\mathbb{Z}_{+}$, при котором
\begin{equation*}
  \frac{f^{(s)}(z_0)}{s!}z_0^s\notin\mathbb{Z}+i\mathbb{Z}.
\end{equation*}

Отметим только, как это следует из формулы Коши-Адамара, что если $\displaystyle f(z)=\sum\limits_{k=0}^\infty c_kz^k$, где все $c_k$ целые, то радиус сходимости ряда, если функция не равна многочлену, не больше единицы.

\textbf{Следствие.} \emph{В обозначениях теоремы 1 при сделанных предположениях}
\begin{equation*}
  E_n^e(f;K_r)=\min\limits_{q_n}\max\limits_{K_r}\big|f(z)-q_n(z)\big|=O\Big(\Big(\frac{r}{R}\Big)^n+(n+1)\rho^n\Big).
\end{equation*}

$\triangleright$ По теореме Бернштейна (см., напр., [\ref{Timan}], [\ref{Dzjad}], [\ref{Trigub_Belinsky}]) существует последовательность $\{p_n\}$ такая, что при $0\leq \nu\leq s$
\begin{equation*}
  \max\limits_{K_r}\Big|f^{(\nu)}(z)-p_n^{(\nu)}(z)\Big|\leq c(r,R,s)\Big(\frac{r}{R}\Big)^n.
\end{equation*}

Так что можно считать, что $p_n^{(\nu)}(z_0)=f^{(\nu)}(z_0)$ $(0\leq\nu\leq s)$.

Отметим ещё, не уменьшая общности, что
\begin{equation*}
  f(z)=\sum\limits_{k=0}^\infty a_k(z-z_0)^k,\qquad \sum\limits_{k=0}^\infty |a_k|\cdot R^k<\infty.
\end{equation*}
Тогда
\begin{equation*}
  \max\limits_{K_r}\Big|\sum\limits_{k=n+1}^\infty a_k(z-z_0)^k\Big|\leq \sum\limits_{k=n+1}^\infty|c_k|\cdot R^k\Big(\frac{r}{R}\Big)^k\leq \Big(\frac{r}{R}\Big)^{n+1} \sum\limits_{k=n+1}^\infty|c_k|\cdot R^k
\end{equation*}
и нужно приблизить, и это главное, $p_n$ многочленом $q_n$.

Имеем
\begin{equation*}
  a_s(z-z_0)^s=a_s(-z_0)^s\Big(\frac{z-z_0}{-z_0}\Big)^s=a_s(-z_0)^s\Big(1-\sum\limits_{k=1}^s b_kz^k\Big)
\end{equation*}
и, как и при доказательстве теоремы 1, заменяем $c_s(-z_0)^s$ на ближайшее целое и повторяем далее доказательство теоремы, повышая степень $z$. \qquad\qquad\qquad\qquad\qquad $\blacktriangleleft$

\begin{center}
    \textbf{\S3 Приближение констант на кубе в $\mathbb{R}^d$}
\end{center}

Пусть $x=(x_1, ..., x_d)\in \mathbb{R}^d$, $\displaystyle x^k=\prod\limits_{j=1}^d x_j^{k_j}$, $k_j\in \mathbb{Z}_{+}$, $|k|=\sum\limits_{j=1}^d k_j$, $p_n(x)=\sum\limits_{0\leq k_j\leq n_j} a_kx^k$, $n=\sum\limits_{j=1}^d n_j$ и
\begin{equation*}
  K=\Pi_{a,b}=\big\{x\in \mathbb{R}^d:\ x_j\in[a,b],\ 1\leq j\leq d\big\}=[a,b]^d.
\end{equation*}

Поскольку куб не должен содержать целочисленных точек, то, не уменьшая общности, считаем $0<a<b<1$ и $O=(0, ..., 0)$ --- ближайшая к $\Pi_{a,b}$ целочисленная точка, т.е. $a+b\leq1$.

Положим при $\lambda\in(0,1)$
\begin{equation*}
  E_n^e(\lambda;\Pi_{a,b})=\min\limits_{q_n}\max\limits_{x\in\Pi_{a,b}}|\lambda-q_n(x)|.
\end{equation*}

\textbf{Теорема 2} \emph{При сделанных предположениях и любом $n\in\mathbb{N}$}
\begin{equation*}
  E_n^e(\lambda;\Pi_{a,b})\leq c(d)n^d\rho^n,\qquad \rho=\max\Bigg\{\frac{\sqrt{b}-\sqrt{a}}{\sqrt{b}+\sqrt{a}},b\Bigg\}
\end{equation*}
\emph{и $\rho$ уменьшить нельзя, в общем случае.}

$\triangleright$ Многочлен Чебышева для отрезка $[a,b]$ равен
\begin{equation*}
  C_n(t;a,b)=C_n(0;a,b)-C_n(0;a,b)\sum\limits_{k=1}^n a_kt^k
\end{equation*}
и при $t\in[a,b]$
\begin{equation*}
  \Big|1-\sum\limits_{k=1}^n a_kt^k\Big|\leq\Big|\frac{C_n(t;a,b)}{C_n(0;a,b)}\Big|\leq\frac{1}{|C_n(0;a,b)|}.
\end{equation*}

Перемножая такие неравенства при $t=x_j$ и $n=n_j$ $(1\leq j\leq d)$, получаем при $x\in\Pi_{a,b}$ в силу \eqref{eq_1}
\begin{equation}\label{eq_2}
  \Big|1-\sum\limits_{1\leq|k|\leq n}b_kx^k\Big|\leq\prod\limits_{j=1}^d\frac{1}{|C_{n_j}(0;a,b)|}\leq2^d\Bigg(\frac{\sqrt{b}-\sqrt{a}}{\sqrt{b}+\sqrt{a}}\Bigg)^n.
\end{equation}

Умножим это неравенство на $\lambda\in(0,1)$ и учтём, что при $\delta_k\in[0,1)$ $\lambda b_k=[\lambda b_k]+\delta_k$. Так что нужно приблизить одночлены $\delta x^s$ при $\delta\in(0,1)$ и $s\in[0,n]$ многочленами $q_n$. Для этого умножим неравенство \eqref{eq_2} на $\delta x^s$, заменяя $n$ на $n-s$.

Получаем
\begin{equation*}
  \Big|\delta x^s-\sum\limits_{s+1\leq|k|\leq n}\delta b_{k-s}x^k\Big|=\Big|\delta x^s-\sum\limits_{1\leq k\leq n-s} \delta b_k x^{k+s}\Big|\leq x^s\cdot2^d\Bigg(\frac{\sqrt{b}-\sqrt{a}}{\sqrt{b}+\sqrt{a}}\Bigg)^{n-s}.
\end{equation*}

Далее выделяем целую часть $[\delta b_{k-s}]$, а к дробной части опять применяем это же неравенство. И так до $s=n$. Приходим к неравенству: при $x\in\Pi_{a,b}$
\begin{equation*}
  |\lambda-q_n(x)|\leq2^d\prod\limits_{s=0}^n x^s\Bigg(\frac{\sqrt{b}-\sqrt{a}}{\sqrt{b}+\sqrt{a}}\Bigg)^{n-s}\leq2^d\prod\limits_{s=0}^n b^s\Bigg(\frac{\sqrt{b}-\sqrt{a}}{\sqrt{b}+\sqrt{a}}\Bigg)^{n-s}\leq2^d\rho^n\textrm{card}\big\{x^k:\ |k|\leq n\big\}.
\end{equation*}

Но, как известно (см., напр., [\ref{Stein_Weiss}], гл. IV, п. 2),
\begin{equation*}
  \textrm{card}\big\{x^k:\ |k|=s\big\}=\binom{d+s-1}{d-1}=\frac{s(s+1)...(s+d-1)}{(d-1)!},
\end{equation*}
откуда
\begin{equation*}
\begin{gathered}
  \textrm{card}\big\{x^k:\ |k|\leq n\big\}=\sum\limits_{s=0}^n\textrm{card}\big\{x^k:\ |k|=s\big\}\leq(n+1)\binom{d+n-1}{d-1}\leq\\
  \leq(n+1)(n+d-1)^{d-1}\leq c_1(d)n^d.
\end{gathered}
\end{equation*}

Оценка приближения сверху доказана.

Для оценки приближения снизу воспользуемся неравенством Бернштейна ([\ref{Bernsht}], II, 93):
\begin{equation*}
  |p_n(0)|\leq\max\limits_{\Pi_{a,b}}|p_n(x)|\prod\limits_{j=1}^d |C_{n_j}(0;a,b)|\leq\max\limits_{\Pi_{a,b}}|p_n(x)|\Bigg(\frac{\sqrt{b}-\sqrt{a}}{\sqrt{b}+\sqrt{a}}\Bigg)^n.
\end{equation*}

Подставляя вместо $p_n$ разность $\lambda-q_n$, получим
\begin{equation*}
  E_n^e(\lambda;\Pi_{a,b})\geq\min\limits_{c\in\mathbb{Z}}|\lambda-c|\Bigg(\frac{\sqrt{b}-\sqrt{a}}{\sqrt{b}+\sqrt{a}}\Bigg)^n
\end{equation*}
и, значит, при $ b\leq a\Big(\frac{1+a}{1-a}\Big)^2$, когда $ b\leq\frac{\sqrt{b}-\sqrt{a}}{\sqrt{b}+\sqrt{a}}$, уменьшить $\rho$ нельзя.

Если же $b> a\Big(\frac{1+a}{1-a}\Big)^2$ и $ b=\frac{1}{q}$ $(q\in\mathbb{N}, q\geq2)$, то
\begin{equation*}
  E_n^e(\lambda;\Pi_{a,b})\geq\min\limits_{q_n}\Big|\lambda-q_n\Big(\frac{p}{q}\Big)\Big|\geq\min\limits_{p_1\in\mathbb{N}}\Big|\lambda-\frac{p_1}{q^n}\Big|\geq\frac{1}{q^{n+2}}=b^{n+2}
\end{equation*}
для бесконечного числа $n$ (см. в \S2).

А при достаточно малом $a$ и $a+b\leq1$.

Теорема 2 доказана. \qquad\qquad\qquad\qquad\qquad\qquad\qquad\qquad\qquad\qquad\qquad\qquad\qquad$\blacktriangleleft$

\vskip 1mm
Случай $\lambda=\frac{p}{q}$ рассмотрен в \S2.

Из теоремы 2 можно вывести следствие, подобное следствию теоремы 1.

\begin{center}
    \textbf{\S4 Интегральные приближения}
\end{center}

Пусть
\begin{equation*}
  K_r=\big\{x=(x_1,...,x_d):\ |x|\leq r\big\}
\end{equation*}
шар радиуса $r$ с центром в нуле, а $n$ --- степень многочлена по совокупности переменных.

\textbf{Теорема 3} \emph{При $\lambda\in(0,1)$, $r\in(0,2)$, $p\in[1,+\infty)$ и $n\rightarrow\infty$}
\begin{equation*}
  E_n^e(\lambda;K_r)=\min\limits_{q_n}\Bigg(\int\limits_{K_r}|\lambda-q_n(x)|^pdx\Bigg)^{\frac{1}{p}}\asymp n^{-\frac{d}{p}}
\end{equation*}
\emph{(двустороннее неравенство с положительными константами, не зависящими от $n$).}

\emph{При $r\geq2$}
\begin{equation*}
  E_n^e(\lambda;K_r)_p\geq c(\lambda,r,p).
\end{equation*}

$\triangleright$ Для оценки приближения снизу понадобится

\textbf{Лемма 1} \emph{При $\alpha>-1$, $p\in[1,+\infty)$ и любой последовательности $\{p_n\}$}
\begin{equation*}
  \int\limits_0^1 t^\alpha|p_n(t)|^p dt\geq c(\alpha,p)|p_n(0)|^p\frac{1}{n^{2+2\alpha}}.
\end{equation*}

$\triangleright$ Если $\displaystyle p_n(t)=\sum\limits_{k=0}^n a_kt^k$, где $a_0=p_n(0)$, то при $p=1$ и $n\geq2$
\begin{equation*}
  \begin{gathered}
  \int\limits_0^1 t^\alpha|p_n(t)|dt\geq\max\limits_{[0,1]}\Big|\int\limits_0^x t^\alpha p_n(t)dt\Big|=\\
  =\max\limits_{[0,1]}x^{\alpha+1}\Big|\sum\limits_{k=0}^n a_k\frac{x^k}{k+\alpha+1}\Big|\geq\frac{1}{n^{2\alpha+2}}\max\limits_{\big[\frac{1}{n^2},1\big]}\Big|\sum\limits_{k=0}^n a_k\frac{x^k}{k+\alpha+1}\Big|.
  \end{gathered}
\end{equation*}

Но в силу экстремального свойства многочленов Чебышева $C_n$ (рост нормы многочлена при расширении отрезка, см. \eqref{eq_1}) при абсолютной константе $c>0$
\begin{equation}\label{eq_3}
  \max\limits_{[0,1]}|p_n(x)|\leq c\max\limits_{\big[\frac{1}{n^2},1\big]}|p_n(x)|.
\end{equation}

Поэтому
\begin{equation*}
  \int\limits_0^1 t^\alpha |p_n(t)|dt\geq \frac{1}{c}\cdot\frac{1}{n^{2+2\alpha}}
\max\limits_{[0,1]}\Big|\sum\limits_{k=0}^n a_k\frac{x^k}{k+\alpha+1}\Big|\geq\frac{1}{c}\cdot\frac{|a_0|}{n^{2+2\alpha}}.
\end{equation*}

Осталось применить неравенство Гёльдера
\begin{equation*}
  \int\limits_0^1 t^\alpha|p_n(t)|dt\leq\Big(\frac{1}{\alpha+1}\Big)^{\frac{1}{p'}}\Big(\int\limits_0^1 t^\alpha|p_n(t)|^pdt\Big)^{\frac{1}{p}}.\qquad\qquad\qquad\qquad\qquad\qquad\qquad \blacktriangleleft
\end{equation*}

Докажем оценку приближения $\lambda$ в теореме снизу.

Обозначим через $S_t$ сферу радиуса $t\in(0,r]$ с центром в нуле $(|x|=t)$, а через $dS_t$ --- элемент её площади. Если при $\displaystyle |k|=\sum\limits_{j=1}^d k_j$ и некоторых целых коэффициентах
\begin{equation*}
  \Big(E_{2n}^e(\lambda;K_r)_p\Big)^p=\int\limits_{|x|\leq r}\Big|\lambda-\sum\limits_{|k|\leq 2n} c_k x^k\Big|^pdx,
\end{equation*}
то, учитывая, что в силу неравенства Гёльдера,
\begin{equation*}
  \Big|\int\limits_{S_t}\Big[\lambda-\sum\limits_{|k|\leq 2n} c_k x^k\Big]dS_t\Big|\leq c(d,p')\Big(\int\limits_{S_t}\Big|\lambda-\sum\limits_{|k|\leq2n} c_kx^k\Big|^pdx\Big)^{\frac{1}{p}},
\end{equation*}
получаем
\begin{equation*}
  \begin{gathered}
  \Big(E_{2n}^e(\lambda;K_r)_p\Big)^p=\int\limits_0^r dt\int\limits_{S_t}\Big|\lambda-\sum\limits_{|k|\leq2n}c_kx^k\Big|^pdS_t
  \geq\frac{1}{c(d,p')}\int\limits_0^r\Big|\int\limits_{S_t}\big[\lambda-c_0-\sum\limits_{1\leq|k|\leq 2n}c_kx^k\big]dS_t\Big|^pdt.
  \end{gathered}
\end{equation*}

Учитывая ещё, что $dS_t=t^{d-1}dS_1$ и симметрию сферы, приходим к неравенству
\begin{equation}\label{eq_4}
  \Big(E_{2n}^e(\lambda;K_r)_p\Big)^p\geq\frac{1}{c(d,p')}\int\limits_0^r t^{d-1}\Big|\lambda-c_0-\sum\limits_{1\leq k\leq n} a_kt^{2k}\Big|^pdt
\end{equation}
и применяем лемму 1 при $ \alpha=\frac{d}{2}-1$ после замены $t^2\rightarrow r^2t$.

Рассмотрим теперь случай $r\geq2$.

В силу неравенства А.~И.~Коркина и Е.~И.~Золотарёва ([\ref{Timan}], [\ref{Trigub_Belinsky}]) при $m\in\mathbb{Z}_{+}$
\begin{equation*}
  \int\limits_a^b t^m\Big|\sum\limits_{k=0}^n a_kt^k\Big|dt\geq 4\Big(\frac{b-a}{4}\Big)^{n+m+1}|a_n|.
\end{equation*}

При $d\geq2$ можно выбрать $\displaystyle\beta\geq-\frac{1}{p}$ так, чтобы
\begin{equation*}
  \Big(\frac{d}{2}-1\Big)^{\frac{1}{p}}+\beta=m\in\mathbb{Z}_{+}.
\end{equation*}

Тогда в силу неравенства Гёльдера для любого многочлена $q_n$
\begin{equation*}
  4\leq\int\limits_0^4 t^m|\lambda-q_n(t)|dt\leq\Big(\int\limits_0^4 t^{\frac{d}{2}-1}|\lambda-q_n(t)|^pdt\Big)^{\frac{1}{p}}\Big(\int\limits_0^4 t^{\beta p'}dt\Big)^{\frac{1}{p'}}
\end{equation*}
и получаем нужную оценку снизу.

При $d=1$ в силу \eqref{eq_4}
\begin{equation*}
  E_{2n}^e(\lambda;K_2)_p\geq\Big(\frac{1}{c(p')}\Big)^{\frac{1}{p}}\Big(\int\limits_{-2}^2\big|\lambda-q_n(t^2)\big|^pdt\Big)^{\frac{1}{p}},
\end{equation*}
а
\begin{equation*}
  \begin{gathered}
  4\leq\int\limits_{-2}^2\big|\lambda-q_n(t^2)\big|dt=\frac{1}{2}\int\limits_0^4\frac{1}{\sqrt{t}}|\lambda-q_n(t)|dt
  \leq\frac{1}{2}\Big(\int\limits_0^4 \frac{1}{\sqrt{t}}|\lambda-q_n(t)|^pdt\Big)^{\frac{1}{p}}\Big(\int\limits_0^4 \frac{dt}{\sqrt{t}}\Big)^{\frac{1}{p'}}.
  \end{gathered}
\end{equation*}

Оценка приближения $\lambda$ снизу (см. ещё \eqref{eq_4}) доказана.

Переходим к оценке приближения числа $\lambda$ сверху.

Пусть $E$ --- компакт в $\mathbb{R}$ с трансфинитным диаметром меньше единицы. Тогда существует многочлен
\begin{equation*}
  P(x)=x^m+..., \qquad \max\limits_{E}|P(x)|<1.
\end{equation*}

Отсюда следует, применяя лемму Больцано-Вейерштрасса, что существует многочлен $X$ с целыми коэффициентами и старшим коэффициентом единица такой, что
\begin{equation*}
  0<\max\limits_E|X(x)|<1.
\end{equation*}

После возведения в квадрат можно считать, что
\begin{equation}\label{eq_5}
  0\leq X(x)\leq\rho<1\quad (x\in E),\quad X(x)=x^l+...\ .
\end{equation}

\textbf{Лемма 2} \emph{Если многочлены $p_n$ делятся без остатка на $X^m$ при некотором $m\in\mathbb{N}$, то можно указать многочлены $q_n$ такие, что}
\begin{equation*}
  \max\limits_E|p_n(x)-q_n(x)|=O\Big(\frac{1}{n^m}\Big).
\end{equation*}

$\triangleright$ Делением $p_n$ можно представить в виде $\Big(N=\Big[\frac{n}{l}\Big]\Big)$
\begin{equation*}
  p_n(x)=\sum\limits_{k=m}^N p_{1,k}(x)X^k(x),
\end{equation*}
где все многочлены $p_{1,k}$ имеют степень не выше $l-1$. После простых преобразований (умножение на единицу, сдвиг в индексах суммирования и перемена порядка суммирования) получаем
\begin{equation*}
  \begin{gathered}
  p_n=\sum\limits_{k=m}^N p_{1,k}X^k\sum\limits_{s=0}^{N-k}\binom{N-k}{s}X^s(1-X)^{N-k-s}=\sum\limits_{k=m}^N\sum\limits_{s=k}^N X^k(1-X)^{N-s}p_{1,k}\binom{N-k}{s-k}=\\
  =\sum\limits_{s=m}^N X^s(1-X)^{N-s}\sum\limits_{k=0}^s p_{1,k}\binom{N-k}{s-k}=\sum\limits_{s=m}^N p_{2,s}X^s(1-X)^{N-s},
  \end{gathered}
\end{equation*}
где $p_{2,s}$ --- многочлены степени не выше $l-1$.

Заменяя коэффициенты $p_{2,s}$ на их целые части, получим многочлены $q_{2,s}$, которые отличаются от $p_{2,s}$ на ограниченную по $n$ величину. Так что
\begin{equation*}
  |p_n(x)-q_n(x)|=O\Big(\sum\limits_{s=m}^N X^s(1-X)^{N-s}\Big).
\end{equation*}

Но при $x\in E$ (см. \eqref{eq_5})
\begin{equation*}
  \sum\limits_{s=m}^N X^s(1-X)^{N-s}\leq\sum\limits_{s=N-m+1}^N X^s+\frac{1}{N^m}\sum\limits_{s=m}^{N-m}\binom{N}{s}X^s(1-X)^{N-s}\leq m\rho^{N-m+1}+\frac{1}{N^m}=O\Big(\frac{1}{n^m}\Big)
\end{equation*}
и лемма 2 доказана.\qquad\qquad\qquad\qquad\qquad\qquad\qquad\qquad\qquad\qquad\qquad\qquad\qquad\qquad $\blacktriangleleft$

\textbf{Лемма 3} \emph{При $\alpha>-1$, $p\in[1,+\infty)$ и $m\in\mathbb{N}$ существует последовательность многочленов $\{p_n\}_m^\infty$ такая, что при некоторой константе $c(\alpha,p,m)$}
\begin{equation*}
  \int\limits_0^1 t^\alpha|p_n(t)|^pdt\leq c(\alpha,p,m)\cdot\frac{1}{n^{2+2\alpha}}
\end{equation*}
\emph{и}
\begin{equation*}
  p_n(0)=1,\quad p_n^{(\nu)}(0)=0\ (1\leq\nu\leq m),\quad \text{а}\ \max\limits_{[0,1]}|p_n(t)|\leq c(\alpha,p,m).
\end{equation*}

$\triangleright$ При $p=2$ соответствующая задача решена давно.

Имеем (см. [\ref{Akhiez}], гл. I, \textbf{14})
\begin{equation*}
  \begin{gathered}
\min\limits_{\{a_k\}_m^n}\int\limits_0^1t^\alpha\Big|1-\sum\limits_{k=m}^n a_kt^k\Big|^2dt=\min\limits_{\{a_k\}_m^n}\int\limits_0^1\Big|t^{\frac{\alpha}{2}}-\sum\limits_{k=m}^n a_kt^{k+\frac{\alpha}{2}}\Big|^2dt=\\
=\frac{1}{\alpha+1}\prod\limits_{k=m}^n\Big(\frac{k}{k+\alpha+1}\Big)^2=\frac{1}{\alpha+1}\sum\limits_{k=m}^n\Big(1+\frac{\alpha+1}{k}\Big)^{-2}.
   \end{gathered}
\end{equation*}

Но при $x\rightarrow0$\ \ $\ln(1+x)=x+O(x^2)$. Поэтому
\begin{equation*}
  \ln\prod\limits_{k=m}^n\Big(1+\frac{\alpha+1}{k}\Big)^2\geq2\Big(\sum\limits_{k=m}^n\frac{\alpha+1}{k}+O(1)\Big)=(2\alpha+2)\ln\frac{n+1}{m}+O(1)
\end{equation*}
и искомое неравенство при $p=2$ доказано.

Выведем из него ограниченность $\{p_n\}_m^n$ в $C[0,1]$.

Учитывая, что при $ t\in\Big[\frac{1}{n^2},1\Big]$ $(n\geq2)$\ \  $\displaystyle t^\alpha\geq\frac{1}{n^{2\alpha}}$, получаем
\begin{equation*}
  \max\limits_{\big[\frac{1}{n^2},1\big]}\int\limits_{\frac{1}{n^2}}^x p_n^2(t)dt\leq c_1(\alpha,m)\cdot\frac{1}{n^2}.
\end{equation*}

Применяя неравенство Маркова для производной многочленов ([\ref{Timan}], [\ref{Trigub_Belinsky}]), имеем
\begin{equation*}
  \max\limits_{\big[\frac{1}{n^2},1\big]}p_n^2(x)\leq c_1(\alpha,m)\frac{(2n+1)^2}{n^2}\cdot\frac{2}{1-\frac{1}{n^2}}\leq c_2(\alpha,m).
\end{equation*}
См. ещё \eqref{eq_3}.

Но эту последовательность $\{p_n\}$ можно использовать при любом $p\in[1,+\infty)$, так как
\begin{equation*}
  \Big(p_{\big[\frac{n}{2}\big]}^2(x)\Big)^p\leq\Big(p_{\big[\frac{n}{2}\big]}(x)\Big)^2\max\limits_{[0,1]}\Big|p_{\big[\frac{n}{2}\big]}(x)\Big|^{2p-2}\leq\big(c\cdot c_2(\alpha,m)\big)^{p-1}\Big(p_{\big[\frac{n}{2}\big]}(x)\Big)^2
\end{equation*}
и все условия при $x=0$ выполнены.\qquad\qquad\qquad\qquad\qquad\qquad\qquad\qquad\qquad\qquad $\blacktriangleleft$

\textbf{Лемма 4} \emph{Пусть $h$ и $m\in\mathbb{N}$, $\{x_\nu\}_1^s\subset[a,b]$, $f\in C^{m+h-1}[a,b]$ и $f^{(k)}(x_\nu)=0$ $(\nu\in[1,s]$, $k\in[0,h-1]$). Тогда существует последовательность $\{p_n\}$ такая, что}
\begin{equation*}
  \int\limits_a^b|f(x)-p_n(x)|^pdx=O\Big(\frac{1}{n^{pm}}+\frac{1}{n^h}\Big)
\end{equation*}
\emph{и}
\begin{equation*}
  p_n^{(k)}(x_\nu)=0\quad (\nu\in[1,s], k\in[0,h]).
\end{equation*}

$\triangleright$ К функции
\begin{equation*}
  F(x)=\frac{f(x)}{\Omega^{d-1}(x)}\in C^m[a,b],\qquad \Omega(x)=\prod\limits_{\nu=1}^s (x-x_\nu)
\end{equation*}
применяем теорему Джексона ([\ref{Timan}], [\ref{Dzjad}], [\ref{Trigub_Belinsky}]): существует последовательность $\{p_n\}$ такая, что
\begin{equation*}
  \int\limits_a^b|F(x)-p_n(x)|^pdx=O\Big(\frac{1}{n^{pm}}\Big),
\end{equation*}
откуда после умножения на $\Omega^{d-1}$
\begin{equation}\label{eq_6}
  \int\limits_a^b|f(x)-\Omega^{d-1}p_n(x)|^pdx=O\Big(\frac{1}{n^{pm}}\Big).
\end{equation}

Считаем далее для простоты, что $a=0$.

В силу леммы 3 существует и ограничена последовательность $\{p_n^{*}\}$ многочленов с условием
\begin{equation*}
  \int\limits_0^b x^{\frac{d}{2}-1}\big|p_n^{*}(x^2)\big|^pdx=O\Big(\frac{1}{n^d}\Big).
\end{equation*}

Поэтому при любом $\nu\in[1,s]$
\begin{equation*}
  \begin{gathered}
\int\limits_0^b |x-x_\nu|^{\frac{d}{2}-1}\big|p_n^{*}\big((x-x_\nu)^2\big)\big|^pdx=\int\limits_{-x_\nu}^{b-x_\nu}|x|^{\frac{d}{2}-1}\big|p_n^{*}(x^2)\big|^pdx\leq\\
\leq\int\limits_{-b}^b |x|^{\frac{d}{2}-1}\big|p_n^{*}(x^2)\big|^pdx=2\int\limits_0^b x^{\frac{d}{2}-1}\big|p_n^{*}(x^2)\big|^pdx=O\Big(\frac{1}{n^d}\Big).
   \end{gathered}
\end{equation*}

Вместо $|x-x_\nu|$ поставим $|\Omega^{d-1}(x)|$ (умножение на ограниченную величину) и многочлен в \eqref{eq_6} изменим на
\begin{equation*}
  \Delta(x)=\sum\limits_{\nu=1}^s p_n(x_\nu)\frac{\Omega(x)}{\Omega'(x_\nu)(x-x_\nu)}p_n^{*}(x-x_\nu),\quad \Delta(x_\nu)=p_n(x_\nu)\ (1\leq\nu\leq s).
\end{equation*}

Так что $p_n(x)-\Delta(x)$ делится без остатка на $\Omega(x)$.

Ещё нужно учесть, что
\begin{equation*}
  \begin{gathered}
\int\limits_0^b\big|\Omega^{d-1}(x)\big|\cdot|\Delta(x)|^pdx=O\Big(\int\limits_0^b|\Omega(x)|^{d-1}\sum\limits_{\nu=1}^s\big|p_n^{*}(x-x_\nu)\big|^pdx\Big)=\\
=O\Big(\sum\limits_{\nu=1}^s\int\limits_0^b|x-x_\nu|^{d-1}\big|p_n^{*}(x-x_\nu)\big|^pdx\Big)=O\Big(\frac{1}{n^d}\Big). \blacktriangleleft
   \end{gathered}
\end{equation*}

Докажем теперь оценку приближения $\lambda$ сверху.

Для отрезка $[0,r^2]$ при $r<2$ существует $X$ с условием \eqref{eq_5}.

Пусть теперь $X_1$ --- многочлен с целыми коэффициентами и простыми нулями $X$. Так что $X_1^l$ делится без остатка на $X$.

При этом
\begin{equation*}
  X_1(x)=X_{1,1}(x)\cdot X_{1,2}(x),
\end{equation*}
где все нули $X_{1,1}$ лежат на $[0,r^2]$, а $X_{1,2}$ не имеет нулей на $[0,r^2]$.

Применяя лемму 4 к функции $\displaystyle f(x)=\lambda\frac{X^m(x)}{X_{1,2}^{m+1}(x)}$ при $m\in \mathbb{Z}_{+}$ и умножая на $X_{1,2}^{m+1}$, получим ($p_{1,k}$ --- многочлены степени $\leq l-1$)
\begin{equation*}
  \int\limits_0^{r^2}\Big|\lambda X^m(x)-\sum\limits_{k=m+1}^{\big[\frac{n}{l}\big]}X^k(x)p_{1,k}(x)\Big|^pdx=O\Big(\frac{1}{n^d}\Big).
\end{equation*}

При $m=0$ заменяем $p_{1,1}$ на многочлен $q_{1,1}$ с целыми коэффициентами так, чтобы $p_{1,1}-q_{1,1}=O(1)$ и применяем это же неравенство
при $m=1$
\begin{equation*}
  \int\limits_0^{r^2}\Big|\lambda -q_{1,1}X(x)-\sum\limits_{k=2}^{\big[\frac{n}{l}\big]}X^k(x)p_{2,k}(x)\Big|^pdx=O\Big(\frac{1}{n^d}\Big).
\end{equation*}

И далее повышаем таким же образом степень $X^m$ до степени $X^d$, чтобы применить лемму 2.
Теорема 3 доказана. \qquad\qquad\qquad\qquad\qquad\qquad\qquad\qquad\qquad\qquad\qquad\qquad$\blacktriangleleft$

\medskip
Отметим, что теорему 3 можно обобщить на тела вращения $(|x|\in E\subset\mathbb{R}_{+})$.

Например, трансфинитный диаметр $E=[-\beta,-\alpha]\cup[\alpha,\beta]$, $0<\alpha<\beta$, равен $ \frac{1}{2}\sqrt{\beta^2-\alpha^2}$ и можно рассмотреть разность двух шаров с центров в нуле
\begin{equation*}
  K_{r,R}=\big\{x\in \mathbb{R}^d,\ 0<r\leq|x|\leq R,\ R^2-r^2<4\big\}.
\end{equation*}

А если$K_{r,R}$ не содержим целочисленных точек, то можно изучить и приближение чисел и гладких функций в $\sup$-норме. А для этого понадобится аналог лемм 1 и 3 при $p=\infty$.

\textbf{Лемма 5} \emph{При $\alpha\geq0$, $m\in\mathbb{N}$ и $n\rightarrow\infty$}
\begin{equation*}
  \min\limits_{p_n}\max\limits_{[0,1]}\Big|t^\alpha\big(1-t^mp_n(t)\big)\Big|\asymp\frac{1}{n^{2\alpha}}
\end{equation*}
\emph{и}
\begin{equation*}
  \min\limits_{p_n}\max\limits_{[-1,1]}|t|^\alpha\Big|1-t^mp_n(t)\Big|\asymp\frac{1}{n^{\alpha}}.
\end{equation*}

$\triangleright$ Случай $\alpha=1$ есть в [\ref{Trigub62}] (выведен из свойств $C_n$). При $\alpha\in\mathbb{N}$ достаточно возвести это соотношение в степень $\alpha$.

Пусть теперь $\alpha>1$ и $\alpha<s$, $s\in\mathbb{N}$.

Для оценки приближения сверху берём $ \varepsilon=\frac{\alpha-1}{s-1}\in(0,1)$. Тогда $\alpha=\alpha_1+\alpha_2$, где $\frac{\alpha_1}{\varepsilon}=s$ и $\frac{\alpha_2}{1-\varepsilon}=1$.

При $\displaystyle A=\max\limits_{[0,1]}\big|1-t^mp_n(t)\big|$ имеем
\begin{equation*}
  t^\alpha A=\Big(t^{\frac{\alpha_1}{\varepsilon}}A\Big)^\varepsilon\Big(t^{\frac{\alpha_2}{1-\varepsilon}}A\Big)^{1-\varepsilon}=O\Big(\frac{1}{n^{2s\varepsilon}}\cdot\frac{1}{n^{2(1-\varepsilon)}}\Big)=O\Big(\frac{1}{n^{2(\alpha_1+\alpha_2)}}\Big).
\end{equation*}
(одновременное возведение в степень).

Пусть теперь $\alpha\in(0,1)$ и $ s=\Big[\frac{1}{\alpha}\Big]+1$.

По доказанному
\begin{equation*}
  \max\limits_{[0,1]}t^{s\alpha}\Big|1-t^mp_n(t)\Big|=O\Big(\frac{1}{n^{2s\alpha}}\Big).
\end{equation*}

На $\displaystyle \Big[\frac{1}{n^2},1\Big]$, а, значит, и на $[0,1]$ (см. \eqref{eq_3})\ \  $1-t^mp_n(t)=O(1)$ и
\begin{equation*}
  \max\limits_{[0,1]}t^{s\alpha}\Big|1-t^mp_n(t)\Big|^s=O\Big(\frac{1}{n^{2s\alpha}}\Big).
\end{equation*}

Остаётся возвести это неравенство в степень $\displaystyle \frac{1}{s}$.

А теперь --- оценка снизу.
\begin{equation*}
  \max\limits_{[0,1]}t^{\alpha}\Big|1-t^mp_n(t)\Big|\geq\frac{1}{n^{2\alpha}}\max\limits_{\big[\frac{1}{n^2},1\big]}\Big|1-t^mp_n(t)\Big|\geq\frac{c(\alpha,m)}{n^{2\alpha}}\max\limits_{[0,1]}\Big|1-t^mp_n(t)\Big|\geq\frac{c(\alpha,m)}{n^{2\alpha}}.
\end{equation*}

Первое соотношение в лемме 5 доказано.

Второе следует из первого при замене $\alpha$ на $\displaystyle\frac{\alpha}{2}$ и $t$ на $t^2$ с учётом чётности.

Лемма 5 доказана. \qquad\qquad\qquad\qquad\qquad\qquad\qquad\qquad\qquad\qquad\qquad\qquad\qquad$\blacktriangleleft$

\textbf{Замечание 1}

В приведенных выше доказательствах (см. \S\S1 и 2) существенную роль при оценке приближения $\lambda$ снизу играют теоремы о росте нормы многочленов $p_n$ вне данного компакта в $\mathbb{C}$ и $\mathbb{R}^d$. Каков максимальный рост по $n$ \ \ $|p_n(x_0)|$ при $x_0\notin K$, если
\begin{equation*}
  \max\limits_{x\in K}|p_n(x)|=1?
\end{equation*}

 В случае комплексной плоскости $\mathbb{C}$ и "хорошего"\ компакта $K$ искомую величину определяют линией уровня функции, конформно и однолистно отображающей внешность $K$ на внешность единичного круга (при специальной нормировке), на которой (линии) находится точка $x_0$ (см., напр., [\ref{Dzjad}], гл. IX).

 Добавим к случаю круга (теорема 1) следующий результат.

 Пусть $r\in(0,1)$ и $ \alpha\in\big(0,\frac{\pi}{2}\big)$, а
 \begin{equation*}
   K_r=K_{r,\alpha}=\big\{x\in\mathbb{C}:\ |z|\leq r,\ \text{Re}~z\leq -r\cos\alpha\big\}.
 \end{equation*}

 Для любого $\lambda\in \mathbb{R}\setminus\mathbb{Z}$ и $\displaystyle r\leq2\sin\frac{\pi\alpha}{2(2\pi-\alpha)}$ при натуральных коэффициентах
 \begin{equation*}
   \lim\limits_{n\rightarrow\infty}\min\limits_{\{c_k\}\in\mathbb{Z}_{+}}\max\limits_{K_r}\Big|\lambda-\sum\limits_{k=0}^n c_kz^k\Big|^{\frac{1}{n}}=2\sin\frac{\pi\alpha}{2(2\pi-\alpha)}.
 \end{equation*}

 Если $\displaystyle r>2\sin\frac{\pi\alpha}{2(2\pi-\alpha)}$ (в частности, $ r=\frac{1}{2}$) и $\lambda$ не является двоично--рациональным числом, тот же предел равен $\displaystyle\frac{1}{2}$.

 \vskip 1mm
 А на прямой $\mathbb{R}$ задача о росте норм давно решена для множеств положительной меры Лебега.

 Е.~Я.~Ремез (1936) доказал, что при $h\in(0,2)$
 \begin{equation*}
   \textrm{meas}\big\{x\in[-1,1]:\ |p_n(x)|\leq1\big\}\leq2-h\quad \Rightarrow\quad \max\limits_{[-1,1]}|p_n(x)|\leq C_n\Big(\frac{2+h}{2-h}\Big),
 \end{equation*}
 где $C_n$ --- многочлен Чебышева для отрезка $[-1,1]$. Равенство имеет место, напр., при $\displaystyle p_n(x)=C_n\Big(\frac{2x+h}{2-h}\Big)$.

 \vskip 1mm
 Отметим ещё, что тогда же, примерно, G.~Polya доказал подобный результат для максимума модуля $p_n^{(n)}$, т.е. старшего коэффициента $p_n$ (см. [\ref{Timan}], \textbf{2.9.13}).

 А недавно [\ref{Tikhonov}] получено точное неравенство о росте для множеств положительной меры на окружности. Там же в [\ref{Tikhonov}] приведен список предшествующих работ о неравенствах типа Ремеза.

 \textbf{Замечание 2} (о многочленах с целыми коэффициентами в анализе).

 \textbf{Первая задача.}

 $K$ --- компакт в $\mathbb{R}$, $\mathbb{C}$ или $\mathbb{R}^d$, а норма в $C$ или $L_p$. Каков $\min\limits_{q_n\not\equiv0}\|q_n\|$?

\vskip 1mm
Ею успешно занимались L.~Kronecker, H.~Minkovskiy, D.~Hilbert, I.~Schur, M.~Fekete, А.~О.~Гельфонд -- Л.~Г.~Шнирельман, Д.~С.~Горшков, E.~Aparisio, автор (см. обзорную статью [\ref{Trigub71}]), F.~Amoroso, Б.~С.~Кашин, P.~Borwein, T.~Erdelyi, I.~E.~Pritsker, автор (см. [\ref{Trigub2019}] со списком литературы).

Этот и следующий список авторов составлен с учётом времени публикации.

\textbf{Вторая задача} (о приближении функций многочленами $q_n$).

Ею успешно занимались, а речь идёт о возможности аппроксимации и порядке приближения в зависимости от функции и степени многочленов, I.~Okada, Р.~О.~Кузьмин, Л.~В.~Канторович, M.~Fekete, G.~Szeg\"{o}, E.~Aparisio, А.~О.~Гельфонд, Г.~А.~Жирнова, E.~Hewitt -- H.~Zuckermen, автор (см. обзор в [\ref{Trigub71}]), L.~B.~O.~Ferquson, С.~Я.~Альпер, M.~Golitschek.

Если функция допускает приближение многочленами $p_n$ на компакте $K\subset\mathbb{C}$, то для приближения многочленами $q_n$ достаточно приблизить константу $ \lambda=\frac{1}{2}$, а для этого необходимо и достаточно существование многочленов $X$ с целыми коэффициентами, удовлетворяющего неравенству
\begin{equation*}
  0<|X(z)|<1\quad (z\in K).
\end{equation*}

В общем случае, если хотя бы одна функция, отличная от многочлена, допускает приближение многочленами $q_n$, то существует многочлен $X$ с условием
\begin{equation*}
  0<\max\limits_K|X(z)|<1
\end{equation*}
и можно ещё считать, что его старший коэффициент равен единице.

Таким образом, получаем первое необходимое условие: трансфинитный диаметр $K$ меньше единицы. А если у таких многочленов $X$ есть обязательные нули на $K$ (целые алгебраические числа вместе с сопряжёнными), то появляются необходимые арифметические условия на функцию (в отличие от интегральной метрики, см. теорему 3 выше). Это к вопросу о возможности аппроксимации (см. [\ref{Lorentz}], гл. 2, \S4).

В работе [\ref{Trigub62}] разработана схема доказательства прямых теорем для гладких функций  (о порядке приближения) на отрезке $[a,b]\subset\mathbb{R}$. Если, напр., это частный случай, $X$ имеет нули только на $K\subset\mathbb{R}$, $0\leq X(x)<1$ на $K$, то функцию приближаем многочленами $p_n$, которые делятся на $X$ (используются арифметические условия на функцию), а затем приближаем многочленами, делящимися на $X^m$ (см. лемму 2 выше в \S4). Это в случае $K\subset\mathbb{R}$.

Но давно получены прямые теоремы о приближении многочленами $p_n$ на компактах в $\mathbb{C}$ (см. [\ref{Dzjad}], гл. IX). А переход к приближениям многочленами $q_n$ сделан пока только для квадрата $[0,1]^2$ (Вит. Волчков). А в случае функций нескольких вещественных переменных, когда вопрос о делимости многочленов существенно усложняется (см. [\ref{Prasolov}]), прямые теоремы получены только для декартовых произведений одномерных компактов [\ref{Trigub_dep}].

Заметим ещё, что в отличие от приближения многочленами $q_n$, в вопросе о приближении целозначными многочленами, которые, по определению, принимают целые значения в целочисленных точках, размеры компакта $K$ не существенны. Так, для компакта, лежащего на $(0,m)$, $m\in\mathbb{N}$, можно взять
\begin{equation*}
  X(x)=\frac{(x-1)(x-2)...(x-m+1)}{(m-1)!}\in(-1,1).
\end{equation*}

Приведём ещё два нерешённых вопроса.

I. По идее Гельфонда-Шнирельмана, используя информацию о наименьших ненулевых нормах $\max\limits_{[0,1]}|q_n(x)|$, не удалось доказать асимптотический закон распределения простых чисел (см., напр., [\ref{Trigub2019}]). Но, возможно, это удастся после перехода к многочленам многих переменных.

II.  В [\ref{Trigub2001}] доказана прямая теорема о приближении многочленами с натуральным коэффициентами на отрезке $[-2,0]$. Как выглядит подобная теорема для $[-3,0]$, напр., когда
\begin{equation*}
  \max\limits_{[-3,0]}\big|x(x+1)(x+2)(x+3)\big(x^2+3x+1\big)\big|=\frac{45}{64}?
\end{equation*}

\end{document}